\documentclass{birkau}

\usepackage[mathscr]{eucal}
\usepackage{enumerate,graphicx,color}
\usepackage{url}
\usepackage{amssymb}
\usepackage{amsmath}
\usepackage{amsthm}
\usepackage{xy}
\usepackage{enumitem}

\numberwithin{equation}{section}

\theoremstyle{plain}
\newtheorem{thm}{Theorem}[section]
\newtheorem{lem}[thm]{Lemma}
\newtheorem{prop}[thm]{Proposition}
\newtheorem{coro}[thm]{Corollary}

\newtheorem{prob}{Problem}

\theoremstyle{definition}
\newtheorem{df}[thm]{Definition}
\newtheorem{rem}[thm]{Remark}
\newtheorem{eg}[thm]{Example}

\usepackage{pgf,tikz}

\usetikzlibrary{calc,positioning,shapes,arrows,arrows.meta}

\tikzset{%
world/.style={circle,draw,minimum size=0.5cm,fill=gray!15},
label/.style={shape=rectangle, inner sep=6pt},
shaded/.style={draw, shape=circle, fill=black!35, inner sep=1.4pt},
unshaded/.style={draw, shape=circle, fill=white, inner sep=1.4pt},
quasi/.style={draw, shape=rectangle, rounded corners=3pt, fill=white, inner sep=2.5pt, minimum height=14.5pt},
blob/.style={draw, shape=rectangle, rounded corners=12pt, thin, densely dotted},
order/.style={thin},
curvy/.style={thin, looseness=1.2, bend angle=70},
fatcurvy/.style={thin, looseness=1.7, bend angle=75},
map/.style={->, densely dashed, shorten >=5pt, shorten <=5pt, >=stealth', looseness=1.1},
operationgj/.style={->, densely dashed, shorten >=5pt, shorten <=18pt, >=stealth', looseness=1.1},
relationlejk/.style={->, shorten >=5pt, shorten <=5pt, >=stealth'},
auto}

\newcommand{\eusbA}{\medsub e {\kern-0.75pt\A\kern-0.75pt}}

\newcommand{\twiddle}[1]{{\smash{\underset{\raise.375ex\hbox{$\smash\sim$}}
       {#1}}\vphantom{\underline{\mathbf{#1}}}}} 
\newcommand{\stwiddle}[1]{\smash{\underset{\smash{\raise.1ex\hbox{\small$\sim$}}}
                         {\mathbf{#1}}}\vphantom{#1}}
\newcommand{\twoT}{\twiddle 2}

\newcommand{\mpe}[2]{\CG^{\rm mp}(#1,#2)}

\newcommand{\cat}[1]{\boldsymbol{\mathscr{#1}}}

\newcommand{\CG}{\cat G}

\newcommand{\neswarrow}{\mathrel{\text{$\nearrow$\llap{$\swarrow$}}}}

\font\bmi=cmmi8 scaled 1440
\newcommand{\powerset}{\raise.6ex\hbox{\bmi\char'175 }}

\DeclareMathOperator{\dom}{dom} 

\renewcommand{\le}{\leqslant}

\renewcommand{\ge}{\geqslant}

\newcommand{\Cld}{\mathrm{Cld}}

\newcommand{\MDFIP}[2]{\langle {\uparrow} #1,{\downarrow} #2\rangle}
\begin{document}
\title{Dual digraphs of finite meet-distributive and modular lattices}

\author[Craig]{Andrew Craig}
\address{Department of Mathematics and Applied Mathematics\\
University of Johannesburg\\PO Box 524, Auckland Park, 2006\\South~Africa\\
and\\
National Institute for Theoretical and Computational Sciences (NITheCS)\\
South Africa}
\email{acraig@uj.ac.za}

\author[Haviar]{Miroslav Haviar}
\address{Department of Mathematics\\
Faculty of Natural Sciences\\
M. Bel University\\
Tajovsk\'{e}ho 40, 974 01 Bansk\'{a} Bystrica\\
Slovakia\\and\\Department of Mathematics and Applied Mathematics\\
University of Johannesburg\\PO Box 524, Auckland Park, 2006\\South~Africa}
\email{miroslav.haviar@umb.sk}

\author[Marais]{Klarise Marais}
\address{Department of Mathematics and Applied Mathematics\\
University of Johannesburg\\PO Box 524, Auckland Park, 2006\\South~Africa}
\email{klarise.marais@gmail.com}





\begin{abstract}
We describe the digraphs that are dual representations of finite lattices satisfying conditions related to
meet-distributivity and
modularity. This is done using the dual digraph representation
of finite lattices
by Craig, Gouveia and Haviar (2015). These digraphs, known as TiRS digraphs, have their origins in the dual representations of lattices by Urquhart (1978) and  Plo\v{s}\v{c}ica (1995). We describe two properties of finite lattices which are weakenings of (upper) semimodularity and lower semimodularity
respectively,
and then show how these properties have a simple description in the dual digraphs. Combined with previous work 
on dual digraphs of semidistributive lattices (2022), 
it 
leads to a dual representation of
finite
meet-distributive
lattices. 
This provides a natural link to finite convex geometries.
In addition, we present
two sufficient conditions
on a finite TiRS
digraph
for its dual lattice to be modular.
We close by posing four open problems.
\end{abstract}

\maketitle


\section{Introduction}\label{sec:intro}

The first dual representation of arbitrary bounded lattices was given by Urquhart in 1978~\cite{U78}. Since then, many different authors have attempted to provide dualities and dual representations of classes of lattices that are not necessarily distributive
(see the recent survey by the first author~\cite{C22}).

In this paper we examine representations for finite lattices that satisfy conditions related to meet-distributivity and modularity. The dual structures
of these finite lattices will be TiRS digraphs satisfying some additional conditions. It was shown by Craig, Gouveia and Haviar~\cite{P3} that there is a one-to-one correspondence between the class of finite lattices and finite digraphs known as TiRS digraphs (see Definition~\ref{def:TiRS} and Theorem~\ref{thm:TiRSrep}). We remark that this correspondence generalises Birkhoff's one-to-one correspondence between finite distributive lattices and finite posets from the 1930s.

We 
introduce and study lattice-theoretic conditions which generalise both lower semimodularity and (upper) semimodularity for finite lattices and seem to be more natural and simpler than the conditions from~\cite{FCA99}. We are also able to provide equivalent conditions to them on the dual TiRS digraph of a finite lattice. We can combine our lattice-theoretic conditions with our previous  
results~\cite{P5}
to characterise the dual digraphs of finite meet-distributive lattices, which correspond to finite convex geometries. 

Currently, the only known dual characterisation of finite modular lattices is given by the theory of Formal Concept Analysis~\cite{FCA99}. A rather complicated condition is available for the standard context dual to a finite semimodular lattice~\cite[Theorem 42]{FCA99}.
We are 
able to provide conditions on the dual digraph of a finite lattice, which are sufficient though not necessary for modularity of the lattice.

The paper is laid out as follows. In Section~\ref{sec:prelim} we provide some background definitions and results that will be needed later on in the paper. Section~\ref{sec:lowsemimod}
defines two conditions which generalise, respectively, (upper) semimodularity and lower semimodularity. We focus on the generalisation of lower semimodularity---a condition we call (JM-LSM) (see Definition~\ref{def:JM-LSM}). We characterise the dual of (JM-LSM) on the dual digraphs of finite lattices. For completeness we state corresponding conditions and results related to upper semimodularity.
In Section~\ref{sec:meet-dist} we combine the results of Section~\ref{sec:lowsemimod} with  results from a recent paper by Craig, Haviar and S\~{a}o Jo\~{a}o~\cite{P5}.
There, characterisations were given of the digraphs dual to finite join- and meet-semidistributive lattices (and hence also finite semidistributive lattices). The combination of these dual characterisations gives us a characterisation of the dual digraphs of finite meet-distributive lattices (also know as locally distributive lattices).
Furthermore, this allows us to describe a new class of structures that is in a one-to-one correspondence with finite convex geometries.
In Section~\ref{sec:mod} we give two sufficient conditions on a finite TiRS digraph for the dual lattice to be modular.
In Section~\ref{sec:conclusion} 
we explicitly list four open problems and
we also indicate why the task of describing digraphs dual to finite modular lattices is challenging.

\section{Preliminaries}\label{sec:prelim}

Central to the representation of a finite lattice that we will use is the notion of a maximal-disjoint filter-ideal pair. This can, equivalently,
be viewed as a maximal partial homomorphism from a lattice $L$ into the two-element lattice.

\begin{df} [{\cite[Section 3]{U78}}]\label{def:DFIP}
Let $L$ be a lattice. 
Then $\langle F, I \rangle$ is a \emph{disjoint filter-ideal pair} of $L$ if $F$ is a filter of $L$ and $I$ is an ideal of $L$ such that
$F \cap I = \varnothing$.
A disjoint filter-ideal pair $\langle F, I \rangle$ is said to be a \emph{maximal disjoint filter-ideal pair} (MDFIP) if there is no disjoint filter-ideal pair $\langle G, J\rangle \neq \langle F, I \rangle$ such that $F \subseteq G$ and $I \subseteq J$.
\end{df}

The following fact was noted by Urquhart. It is needed for our characterisation of MDFIPs in Theorem~\ref{thm:charMDFIPs}.
\begin{prop}[{\cite[p.~52]{U78}}]\label{prop:MDFIP-JM}
Let $L$ be a finite lattice. If $\langle F, I \rangle$ is an MDFIP 
of $L$ then $\bigwedge F$ is join-irreducible
and $\bigvee I$ is meet-irreducible.
\end{prop}

The set of
join-irreducible
elements of $L$ is denoted $\mathsf{J}(L)$ and the set of
meet-irreducible elements is denoted $\mathsf{M}(L)$.

Given a lattice $L$, we will add a set of arcs to the set of MDFIPs of $L$. The use of such digraphs for lattice representation is due to Plo\v{s}\v{c}ica~\cite{Plos95}. We point out that the original work using (topologised) digraphs used so-called \emph{maximal partial homomorphisms} (see~\cite[Section 1]{Plos95}). It is easy to show that these are in a one-to-one correspondence with MDFIPs.
%

For a lattice $L$, we now present its  dual digraph $G_L = ( X_L, E)$ where the vertices are the  MDFIPs of $L$.
Plo\v{s}\v{c}ica's relation $E$, when transferred to to the set of MDFIPs, is defined  below for two MDFIPs $\langle F, I \rangle$ and $\langle G, J \rangle$:
\begin{itemize}
\item[(E)] $\langle F, I \rangle E \langle G, J \rangle \quad
\iff \quad F \cap J = \emptyset$.
\end{itemize}

For finite lattices every filter is the up-set of a unique element and every ideal is the down-set of a unique element, so we can represent every disjoint filter-ideal pair $\langle F, I \rangle$ by an
ordered pair $\langle {\uparrow}a, {\downarrow}b\rangle$ where $a = \bigwedge F$ and $b= \bigvee I$. Hence for finite lattices we have $\MDFIP{a}{b} E \MDFIP{c}{d}$ if and only if $a \nleqslant d$.
For a digraph $G=(V,E)$ we
let
$xE = \{\, y \in V \mid
xEy
\,\}$ and $Ex = \{\, y \in V \mid
yEx
\,\}$.
The next lemma is easy to prove and it will be useful later on.
\begin{lem}\label{lem:xEyE}
Let $G_L=(X_L,E)$ be the dual digraph of a finite lattice $L$. If $x=\MDFIP{a}{b}$ and $y=\MDFIP{c}{d}$, then
\begin{enumerate}[label=\normalfont(\roman*)]
\item $xE \subseteq yE$ if and only if $a \leqslant c$;
\item $Ex \subseteq Ey$ if and only if $d \leqslant b$.
\end{enumerate}
\end{lem}
%

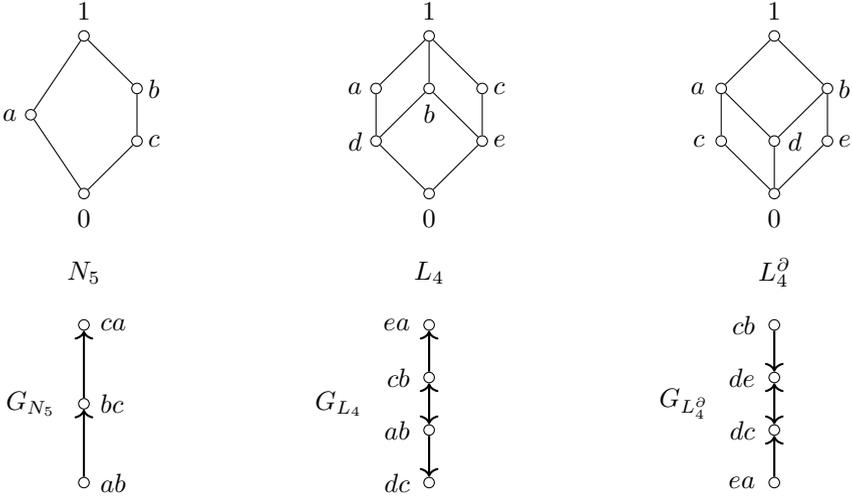
\begin{figure}[ht]
\centering
\begin{tikzpicture}[scale=0.7]
\begin{scope}[xshift=-8.5cm]
\node[unshaded] (bot) at (0,0) {};
\node[unshaded] (a) at (-1,1.5) {};
\node[unshaded] (b) at (1,2) {};
\node[unshaded] (c) at (1,1) {};
\node[unshaded] (top) at (0,3) {};
\draw[order] (bot)--(a)--(top)--(b)--(c)--(bot);
\node[label,anchor=north] at (bot) {$0$};
\node[label,anchor=east,xshift=1pt] at (a) {$a$};
\node[label,anchor=west,xshift=-2pt] at (b){$b$};
\node[label,anchor=west,xshift=-2pt] at (c){$c$};
\node[label,anchor=south] at (top) {$1$};
\node at (0,-1.5) {$N_5$};
\end{scope}
\begin{scope}[xshift=-2cm]
\node[unshaded] (bot) at (0,0) {};
\node[unshaded] (a) at (-1,2) {};
\node[unshaded] (b) at (0,2) {};
\node[unshaded] (c) at (1,2) {};
\node[unshaded] (d) at (-1,1) {};
\node[unshaded] (e) at (1,1) {};
\node[unshaded] (top) at (0,3) {};
\draw[order] (top)--(a)--(d)--(bot);
\draw[order] (b)--(e);
\draw[order] (d)--(b)--(top)--(c)--(e)--(bot);
\node[label,anchor=east,xshift=1pt] at (a) {$a$};
\node[label,anchor=north] at (b) {$b$};
\node[label,anchor=west,xshift=-2pt] at (c) {$c$};
\node[label,anchor=east,xshift=1pt] at (d) {$d$};
\node[label,anchor=west,xshift=-2pt] at (e) {$e$};
\node[label,anchor=north] at (bot) {$0$};
\node[label,anchor=south] at (top) {$1$};
\node at (0,-1.5) {$L_4$};
\end{scope}

\begin{scope}[xshift=4.5cm]
\node[unshaded] (bot) at (0,0) {};
\node[unshaded] (a) at (-1,2) {};
\node[unshaded] (b) at (1,2) {};
\node[unshaded] (c) at (-1,1) {};
\node[unshaded] (d) at (0,1) {};
\node[unshaded] (e) at (1,1) {};
\node[unshaded] (top) at (0,3) {};
\draw[order] (bot)--(c)--(a)--(d)--(b)--(top)--(a);
\draw[order] (d)--(bot)--(e)--(b);
\node[label,anchor=north] at (bot) {$0$};
\node[label,anchor=east] at (a) {$a$};
\node[label,anchor=west,xshift=-2pt] at (b) {$b$};
\node[label,anchor=east] at (c) {$c$};
\node[label,anchor=west,xshift=-1pt] at (d) {$d$};
\node[label,anchor=west,xshift=-2pt] at (e) {$e$};
\node[label,anchor=south] at (top) {$1$};
\node at (0,-1.5) {$L_4^\partial$};
\end{scope}


\begin{scope}[yshift=-5.5cm,xshift=-8.5cm]
\node[unshaded] (ab) at (0,0) {};
\node[label,anchor=west] at (ab) {$ab$};
\node[unshaded] (bc) at (0,1.5) {};
\node[label,anchor=west] at (bc) {$bc$};
\node[unshaded] (ca) at (0,3) {};
\node[label,anchor=west] at (ca) {$ca$};
\path[thick,->] (ab.north) edge  (bc.south);
\path[thick,->] (bc.north) edge  (ca.south);
\node at (-1,1.5) {$G_{N_5}$};
\end{scope}

\begin{scope}[yshift=-5.5cm,xshift=-2cm]
\node[unshaded] (dc) at (0,0) {};
\node[unshaded] (ab) at (0,1) {};
\node[unshaded] (cb) at (0,2) {};
\node[unshaded] (ea) at (0,3) {};
\node[label,anchor=east,xshift=-1pt] at (dc) {$dc$};
\node[label,anchor=east,xshift=-1pt] at (ab) {$ab$};
\node[label,anchor=east,xshift=-1pt] at (cb) {$cb$};
\node[label,anchor=east,xshift=-1pt] at (ea) {$ea$};
\path[thick,<-] (dc.north) edge  (ab.south);
\path[thick,<->] (ab.north) edge  (cb.south);
\path[thick,->] (cb.north) edge (ea.south);
\node at (-1.7,1.5) {$G_{L_4}$};
\end{scope}

\begin{scope}[yshift=-5.5cm,xshift=4.5cm]
\node[unshaded] (ea) at (0,0) {};
\node[unshaded] (dc) at (0,1) {};
\node[unshaded] (de) at (0,2) {};
\node[unshaded] (cb) at (0,3) {};
\node[label,anchor=east, xshift=-1pt] at (ea) {$ea$};
\node[label,anchor=east,xshift=-1pt] at (dc) {$dc$};
\node[label,anchor=east,xshift=-1pt] at (de) {$de$};
\node[label,anchor=east,xshift=-1pt] at (cb) {$cb$};
\path[thick,->] (ea.north) edge  (dc.south);
\path[thick,<->] (dc.north) edge  (de.south);
\path[thick,->] (cb.south) edge  (de.north);
\node at (-1.7,1.5) {$G_{L_4^{\partial}}$};
\end{scope}

\end{tikzpicture}
\caption{Finite lattices $N_5$, $L_4$, $L_4^{\partial}$ and their dual digraphs.}\label{fig:dual-examples}
\end{figure}

Figure~\ref{fig:dual-examples} shows three lattices and their dual digraphs. These three examples will be important throughout this paper.  To make the labelling more succinct, we have denoted by $ab$ the MDFIP $\MDFIP{a}{b}$.
We have also left out the loop on each vertex to keep the display less cluttered. Observe that  the directed edge set is not a transitive relation. The labels $L_4$ and $L_4^\partial$ (as well as $L_3^{\partial}$ which appears later) come from the paper by Davey et al.~\cite{DPR-75}.

The digraphs coming from lattices were described by Craig, Gouveia and Haviar~\cite{P3}.

\begin{df}[{\cite[Definition 2.2]{P3}}]\label{def:TiRS}
A TiRS digraph
$G = (  V,E)$ is a set $V$ and a reflexive
relation $E \subseteq V\times V$ such that:
\begin{itemize}
\item[(S)] If $x,y \in V$ and $x\neq y$ then $xE\neq yE$ or $Ex \neq Ey$.
\item[(R)] For all $x,y \in V$, if $xE \subset yE$ then $(x,y)\notin E$, and if $Ey \subset Ex$ then $(x,y)\notin E$.
\item[(Ti)] For all $x,y \in V$, if $xEy$ then there exists $z \in V$ such that $zE \subseteq xE$ and $Ez \subseteq Ey$.
\end{itemize}
\end{df}

The result below gives a description of dual digraphs of lattices.
\begin{prop}[{\cite[Proposition 2.3]{P3}}]
\label{prop:dualgraphisTiRS}
For any bounded lattice $L$, its dual digraph $G_L=(X_L, E)$ is a TiRS digraph.
\end{prop}

We recall from~\cite{Plos95} a fact concerning general graphs $G =(X,E)$. Let $\twoT=(\{0,1\}, \le)$ denote the two-element graph. A~partial map $\varphi \colon X \to \twoT$ preserves the relation $E$ if  $x, y \in \dom \varphi$ and 
$xEy$
imply $\varphi(x)\leqslant \varphi(y)$. The
set of
maximal  partial $E$-preserving maps (i.e. those that cannot be properly extended) from
$G$
to $\twoT$ is denoted by
$\mpe{G}{\twoT}$.
We use the abbreviation MPEs for 
such partial maps. 

For a graph $G=(X,E)$
and $\varphi, \psi \in \mpe{G}{\twoT}$,
it was shown by Plo\v{s}\v{c}ica~\cite[Lemma 1.3]{Plos95} that
$\varphi^{-1}(1) \subseteq \psi^{-1}(1) \:\Longleftrightarrow\: \psi^{-1}(0) \subseteq \varphi^{-1}(0)$.
This implies that the reflexive and transitive binary relation $\leqslant$ defined on $\mpe{G}{\twoT}$ by
$\varphi\le\psi \iff \varphi^{-1}(1)\subseteq\psi^{-1}(1)$
is a partial order.
In fact, this is a lattice order~\cite[Theorem 2.3]{P1}.
For a graph $G=(X,E)$, denote by $\mathbb{C}(G)$ the (complete) lattice of MPEs $(\mpe{X}{\twoT},\subseteq)$.

The theorem below gives a one-to-one correspondence between finite lattices and finite TiRS digraphs.
This result is essential to the work done in the rest of the current paper.
\begin{thm}[{\cite[Theorem 1.7 and p.~87]{P3}}] \label{thm:TiRSrep}
For any finite bounded lattice $L$ we have that $L$ is isomorphic to $\mathbb{C}(G_L)$ and for any finite
TiRS digraph  $G=(V,E)$
we have that $G$ is isomorphic to  $G_{\mathbb{C}(G)}$.
\end{thm}

\section{Generalising lower and upper semimodularity}
\label{sec:lowsemimod}

For lattice elements $a$ and $b$ we write
$a \prec b$ to denote that $a$ is covered
by~$b$. A lattice is \emph{upper semimodular} if whenever $a \wedge b \prec a$ then $b \prec a \vee b$. It is common to refer to such lattices as \emph{semimodular}. A lattice is \emph{lower semimodular} if whenever $a \prec a \vee b$ then $a \wedge b \prec b$. We use (USM) and (LSM) as abbreviations for these two conditions.

The lattices in Figure~\ref{fig:dual-examples} provide useful examples: $N_5$ satisfies neither (USM) nor (LSM), $L_4$ satisfies
(USM) but not (LSM),
and $L_4^\partial$ satisfies
(LSM) but not (USM).

We will focus on lower semimodularity, rather than upper semimodularity, because of the connection between lower semimodularity and finite convex geometries (see Section~\ref{sec:meet-dist}).
We note that modularity implies both semimodularity and lower semimodularity. If a lattice $L$ has finite length and is semimodular and lower semimodular, then $L$ is also modular (cf.~\cite[Corollary 376]{LTF}).
For further reading we refer to the book by Stern~\cite{S99}.

\begin{figure}[ht]
\centering
\begin{tikzpicture}[scale=1.2]
\begin{scope}[yshift=-5cm,xshift=4cm]
\node[unshaded] (m) at (-1,0) {};
\node[unshaded] (a) at (-1,2) {};
\node[unshaded] (c) at (1,2) {};
\node[unshaded] (b) at (0,2) {};
\node[unshaded] (e) at (0,1) {};
\node[unshaded] (top) at (0,3) {};
\node[unshaded] (d) at (-2, 1) {};
\node[unshaded] (f) at (2, 1) {};
\node[unshaded] (g) at (1,0) {};
\node[unshaded] (h) at (-1,1) {};
\node[unshaded] (i) at (1,1) {};
\node[unshaded] (j) at (0,0) {};
\node[unshaded] (k) at (-3,0) {};
\node[unshaded] (l) at (3,0) {};
\node[unshaded] (bot) at (0,-3) {};
\node[unshaded] (n) at (0,-1) {};
\node[unshaded] (o) at (-1,-1) {};

\draw[order] (b)--(top)--(c)--(e)--(a)--(top);
\draw[order] (e)--(m)--(d)--(a);
\draw[order] (e)--(g)--(f)--(c);
\draw[order] (b)--(i)--(j)--(h)--(b);
\draw[order] (d)--(k)--(bot)--(l)--(f);
\draw[order] (k)--(h);
\draw[order] (l)--(i);
\draw[order] (g)--(n)--(m);
\draw[order] (j)--(o)--(bot)--(n);

\node[label,anchor=north] at (bot) {\small distributive};
\node[label,anchor=west,xshift=0pt] at (n) {\small modular};
\node[label,anchor=east,xshift=1pt] at (o) {B};
\node[label,anchor=west,xshift=0pt] at (j) {\small SD};
\node[label,anchor=east,xshift=0.2cm,yshift=0.2cm] at (h) {\small JSD};
\node[label,anchor=west,xshift=-0.2cm,yshift=0.2cm] at (i) {\small MSD};
\node[label,anchor=east,xshift=1pt] at (a) {\small $g \nearrow m, g \swarrow n \Rightarrow g \neswarrow m$ };
\node[label,anchor=west,xshift=-1pt] at (c) {\small $g \swarrow m, h \nearrow m \Rightarrow g \neswarrow m$};
\node[label,anchor=south,xshift=0pt] at (b) {\small semi-convex};
\node[label,anchor=east,xshift=1pt] at (d) {\small LSM};
\node[label,anchor=west,xshift=1pt] at (f) {\small USM};
\node[label,anchor=south,xshift=-0.7cm,yshift=-0.5cm] at (k) {\small meet-};
\node[label,anchor=south,xshift=-0.3cm,yshift=-0.9cm] at (k) {\small distributive};
\node[label,anchor=south,xshift=0.3cm,yshift=-0.7cm] at (l) {\small join-};
\node[label,anchor=south,xshift=0.75cm,yshift=-1.05cm] at (l) {\small distributive};
\end{scope}

\end{tikzpicture}
\caption{Relationships between generalisations of distributivity.}
\label{fig:FCAbookdiagram}
\end{figure}
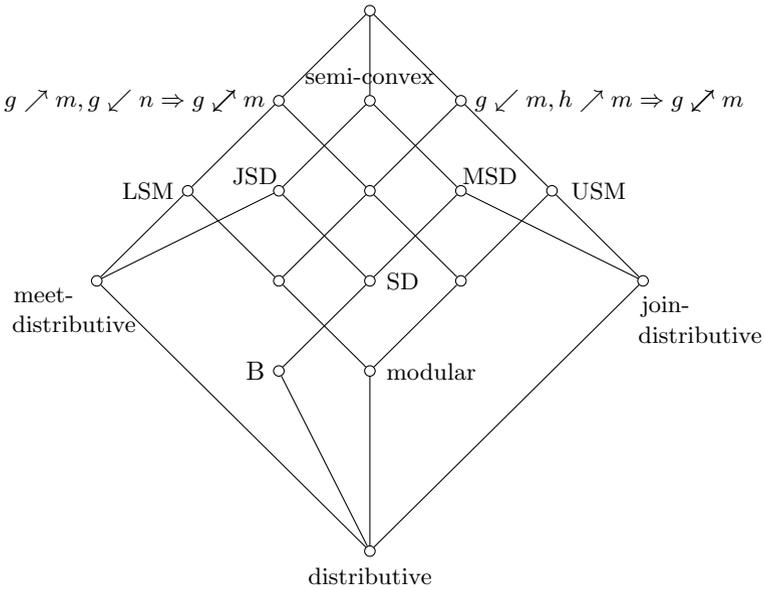

Figure~\ref{fig:FCAbookdiagram} presents a number of different generalisations of distributivity and modularity (including those presented above) and the relationships between them. Observe that the conditions in the top left and top right, which are weakenings of (LSM) and (USM) respectively, are in fact conditions on the standard context dual to a finite lattice. For the necessary terms and notation, we  refer to the book
from where Figure~\ref{fig:FCAbookdiagram} is taken~\cite[p.~234]{FCA99}.

We begin by proving some new results about MDFIPs. 
These will be needed in the proofs of later results.

\begin{lem}\label{lem:covers-imply-maximal}
Let $L$ be a finite lattice.
\begin{enumerate}[label=\normalfont(\roman*)]
\item If $b \in \mathsf{M}(L)$ and $b \prec a \vee b$, then ${\downarrow}b$ is maximal with respect to being disjoint from ${\uparrow}a$.
\item If $a \in \mathsf{J}(L)$ and $a \wedge b \prec a$, then ${\uparrow}a$ is maximal with respect to being disjoint from ${\downarrow}b$.
\end{enumerate}
\end{lem}
\begin{proof}
Assume that $b \in \mathsf{M}(L)$ and $b \prec a \vee b$. This implies $b < a\vee b$ and hence $ a \nleqslant b$ and so ${\uparrow}a \cap {\downarrow}b = \emptyset$. Suppose the ideal ${\downarrow}b$ were to be extended to ${\downarrow}c$ with $b<c$
and ${\uparrow}a \cap {\downarrow}c = \emptyset$.
Since  $b \in \mathsf{M}(L)$, the element $a\vee b$ is the unique upper cover of $b$ and so $a\vee b \in {\downarrow}c$. This implies
$a\vee b \in {\uparrow}a \cap {\downarrow}c$, a contradiction,
showing the maximality of ${\downarrow} b$
with respect to being disjoint from ${\uparrow}a$.

Now assume that $a \in \mathsf{J}(L)$ and $a \wedge b \prec a$. Since $a \wedge b < a$ we have $a \nleqslant b$ and so ${\uparrow}a \cap {\downarrow}b= \emptyset$. If ${\uparrow} a$ were extended to ${\uparrow}d$ with $d < a$
and ${\uparrow}d \cap {\downarrow}b = \emptyset$,
then $d \leqslant a \wedge b$ (the unique lower cover of $a$). We get $a \wedge b \in {\uparrow}d\cap {\downarrow}b$, which shows that ${\uparrow}a$ is maximal with respect to being disjoint from ${\downarrow}b$.
\end{proof}

The next
theorem
gives a characterisation of MDFIPs.

\begin{thm}\label{thm:charMDFIPs}
A disjoint filter-ideal pair $\MDFIP{a}{b}$ is an MDFIP if and only if
it satisfies the following conditions:
\begin{enumerate}[label=\normalfont(\roman*)]
\item $a \in \mathsf{J}(L)$;
\item $b \in \mathsf{M}(L)$;
\item $b \prec a \vee b$;
\item $a \wedge b \prec a$.
\end{enumerate}
\end{thm}
\begin{proof} If $\MDFIP{a}{b}$ is an MDFIP, by Proposition~\ref{prop:MDFIP-JM}, $a \in \mathsf{J}(L)$ and $b \in \mathsf{M}(L)$.
We also have  $b <a\vee b$, since $b=a\vee b$ would imply $a \in {\downarrow}b$. Suppose there exists $c\in L$ such that $b < c <a \vee b$. If $a \leqslant c$ then $c$ would be an upper bound for $\{a,b\}$ and then $a \vee b \leqslant c$. Therefore $a \nleqslant c$.  This would make
$\MDFIP{a}{c}$ a disjoint filter-ideal pair with ${\downarrow}b \subsetneq {\downarrow}c$, contradicting the maximality of the pair
$\MDFIP{a}{b}$.
A dual argument can be applied to show that $a \wedge b \prec a$.

Assume $\MDFIP{a}{b}$ satisfies (i)--(iv).
Lemma~\ref{lem:covers-imply-maximal} then tells us that ${\downarrow}b$ is maximal with respect to being disjoint from ${\uparrow}a$ and vice versa. Hence $\MDFIP{a}{b}$ is an MDFIP.
\end{proof}

The lemmas below will be used in our later investigations.

\begin{lem}\label{lem:JM}
Let $L$ be a finite lattice, $a,b \in L$.
Then the following are equivalent:
\begin{enumerate}[label=\normalfont(\roman*)]
\item $a \nleqslant b$;
\item there exists $j \in \mathsf{J}(L)$ such that $j \leqslant a$ and $j \nleqslant b$;
\item there exists $m \in \mathsf{M}(L)$ such that $b \leqslant m$ and $a \nleqslant m$.
\end{enumerate}
\end{lem}

\begin{proof}
It is well-known that in a finite lattice the set $\mathsf{J}(L)$ is join-dense.
Hence
$a \le b$ is equivalent to the condition
that for all $j \in \mathsf{J}(L)$, $j \le a$ implies $j \le b$.
This settles the equivalence of (i) and (ii).
The equivalence of (i) and (iii) follows similarly from the meet-density of $\mathsf{M}(L)$ in $L$.
\end{proof}

For $a,b \in L$ we define the set $T_{ab}:=\{\, m \in \mathsf{M}(L) \mid b \leqslant m, a \nleqslant m\,\}$. An important consequence of Lemma~\ref{lem:JM} is that $T_{ab}$ is non-empty whenever $a \nleqslant b$. This is needed for our next result.
\begin{lem}\label{lem:Tab}
Let $L$ be a finite lattice and $a, b \in L$, $a \nleqslant b$. Let
$d$ be a maximal element of $T_{ab}$. Then
$d \prec d \vee a$.
\end{lem}

\begin{proof}
Firstly, we point out that $T_{ab}$ is a non-empty finite poset and hence has a maximal element.
Since $a \nleqslant d$, we have $a \vee d \neq d$, and so $d < d \vee a$. Suppose there exists $c\in L$ such that $d < c < d \vee a$. As $d \vee a \nleqslant c$,
by Lemma~\ref{lem:JM}
there exists $m \in \mathsf{M}(L)$ such that
$c \leqslant m$ but $d \vee a \nleqslant m$.
So $d<m$.
If $a \leqslant m$ then
$d \vee a \leqslant m$. It follows $a \nleqslant m$ and $b \le d < m$, so $m \in T_{ab}$. Since $d$ was maximal in $T_{ab}$ and $d < m$, we get a contradiction. Hence $d \prec d \vee a$.
\end{proof}

From the previous lemmas one can derive the following result.
\begin{prop}\label{prop:ambj} Let $L$ be a finite lattice with $a \in \mathsf{J}(L)$ and $b\in\mathsf{M}(L)$. Then
\begin{enumerate}[label=\normalfont(\roman*)]
\item there exists $m \in \mathsf{M}(L)$ such that $\MDFIP{a}{m}$ is an MDFIP;
\item there exists $j \in \mathsf{J}(L)$ such that $\MDFIP{j}{b}$ is an MDFIP.
\end{enumerate}
\end{prop}
\begin{proof} We prove only (i), as then (ii) will follow by a dual argument. Since $a \in \mathsf{J}(L)$, it has a unique lower cover $c$. Clearly $a \nleqslant c$, so by Lemma~\ref{lem:Tab}, there exists a maximal element $m \in T_{ac}$ such that $m \prec m \vee a$. From Lemma~\ref{lem:covers-imply-maximal}(i) we know that ${\downarrow}m$ is maximal with respect to being disjoint from ${\uparrow}a$. If it
were possible
to extend ${\uparrow}a$ to ${\uparrow}d$ with $d <a$, then since $c$ is the unique lower cover of $a$, we would get $c \in {\uparrow}d \cap {\downarrow}m$. Hence ${\uparrow}a$ is maximal with respect to being disjoint from ${\downarrow}m$. It follows that  $\MDFIP{a}{m}$ is an MDFIP.
\end{proof}

We now define a new condition, (JM-LSM), which will be central to the results that follow.
We believe it is a more natural weakening of (LSM) than the condition given in the top left of Figure~\ref{fig:FCAbookdiagram}.
The name of the condition comes from the fact that it is almost identical to the condition (LSM), but the elements involved are quantified over $\mathsf{J}(L)$ and $\mathsf{M}(L)$.

\begin{df}\label{def:JM-LSM}
A finite lattice $L$ satisfies (JM-LSM) if for any $a \in \mathsf{J}(L)$ and $b \in \mathsf{M}(L)$, if $b \prec a \vee b$ then $a \wedge b \prec a$.
\end{df}

\begin{figure}[ht]
\centering
\begin{tikzpicture}[scale=0.7]
\begin{scope}[yshift=-2.5cm,xshift=4cm]
\node[unshaded] (top) at (0,0) {};
\node[unshaded] (a) at (-1,-2) {};
\node[unshaded] (c) at (1,-2) {};
\node[unshaded] (b) at (0,-2) {};
\node[unshaded] (d) at (0,-1) {};
\node[unshaded] (bot) at (0,-3) {};
\draw[order] (d)--(b)--(bot)--(c)--(top)--(d)--(a)--(bot);
\node[label,anchor=north] at (bot) {$0$};
\node[label,anchor=east,xshift=1pt] at (a) {$a$};
\node[label,anchor=west,xshift=-1pt] at (c) {$c$};
\node[label,anchor=west,xshift=-2pt] at (b){$b$};
\node[label,anchor=east,xshift=1pt] at (d){$d$};
\node[label,anchor=south] at (top) {$1$};
\end{scope}

\begin{scope}[yshift=-6cm,xshift=10cm]
\node[unshaded] (ac) at (0,0.5) {};
\node[label,anchor=north] at (ac) {$ac$};
\node[unshaded] (bc) at (2,0.5) {};
\node[label,anchor=north] at (bc) {$bc$};
\node[unshaded] (ab) at (-0.5,2.5) {};
\node[label,anchor=east] at (ab) {$ab$};
\node[unshaded] (ba) at (2.5,2.5) {};
\node[label,anchor=west] at (ba) {$ba$};
\node[unshaded] (cd) at (1,3.5) {};
\node[label,anchor=south] at (cd) {$cd$};
\path[thick,shorten <=2.5pt, shorten >=2.5pt,<->] (ac.east) edge  (bc.west);
\path[thick,<->,shorten >=2pt,shorten <=2pt] (ac.north) edge (ab.south);
\path[thick,<->,shorten >=2pt,shorten <=2pt] (bc.north) edge (ba.south);
\path[thick, shorten >=3pt, shorten <=4pt, ->] (ab.east) edge (bc.north west);
\path[thick,shorten >=3pt, shorten <=4pt,->] (ba.west) edge  (ac.north east);
\path[thick, shorten >=3.5pt, shorten <=3pt, <-] (ab.north east) edge (cd.south west);
\path[thick, shorten >=3.5pt, shorten <=3pt, <-] (ba.north west) edge (cd.south east);
\end{scope}

\end{tikzpicture}
\caption{A finite lattice that satisfies
(JM-LSM) but not (LSM).
Its dual digraph (right) satisfies (LTi).
}\label{fig:L3dualdigraph}
\end{figure}
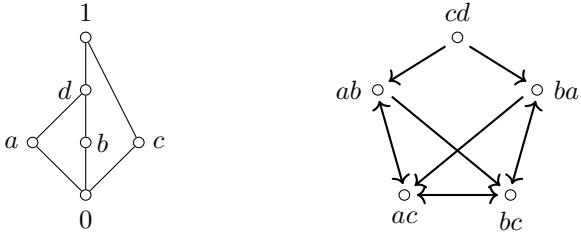

\begin{eg}
Condition (JM-LSM)
is a proper weakening of the condition (LSM).
Indeed,
the lattice in Figure~\ref{fig:L3dualdigraph} satisfies (JM-LSM) but not (LSM).
To see this, observe that
$c \prec c \vee d$ and $c \wedge d \nprec d$, yet $d\notin \mathsf{J}(L)$.

We note that the lattice $L_4$ in Figure~\ref{fig:dual-examples}
does not satisfy (LSM), and
also does not satisfy (JM-LSM):
$c \in \mathsf{J}(L)$, $a \in \mathsf{M}(L)$ and $a \prec c \vee a$, yet $c \wedge a \nprec c $.
\end{eg}

Below
is a condition that we will prove is
equivalent to (JM-LSM).
It
will assist us in proving that
the digraph condition (LTi), given in Definition~\ref{def:LTi},
can be used to characterise
the dual digraphs of finite 
(JM-LSM) lattices.

\begin{df}\label{def:labc}
Condition (L-abc): Let $a \in \mathsf{J}(L)$ and $b \in \mathsf{M}(L)$. If $a \nleqslant b$ then there exists $c \ge b$ such that $\MDFIP{a}{c}$ is an MDFIP.
\end{df}
Notice that if $\MDFIP{a}{c}$ is an MDFIP, then Proposition~\ref{prop:MDFIP-JM}
(cf. also Theorem~\ref{thm:charMDFIPs})
implies that for the element $c$ in Definition~\ref{def:labc} we have $c \in \mathsf{M}(L)$.
Notice also that the
finite lattice
$L_4$ in Figure~\ref{fig:dual-examples}
does not satisfy (L-abc):
we have $a \in \mathsf{J}(L)$, $c \in \mathsf{M}(L)$ and $a \nleqslant c$ and there is no $m \ge c$ such that $\MDFIP{a}{m}$ is an MDFIP.

The following theorem shows that for finite lattices the central property {\upshape (JM-LSM)} can be characterised exactly via the condition {\upshape (L-abc)}.

\begin{thm}\label{thm:labc-iff-JMLSM}
A finite lattice satisfies
{\upshape (JM-LSM)} iff it satisfies {\upshape (L-abc)}.
\end{thm}

\begin{proof}
Assume
(JM-LSM) and let $a \in \mathsf{J}(L)$, $b \in \mathsf{M}(L)$ and $a \nleqslant b$. Let 
$T_{ab}= \{m \in \mathsf{M}(L) \mid b \le m \And a\nleqslant m \}$. 
Then $T_{ab}$ is
a non-empty finite poset.
Hence it has a maximal element, say $c$. So
$c  \in \mathsf{M}(L)$, $b\le c$ and $\MDFIP{a}{c}$ is a disjoint filter-ideal pair. 
To show that $\MDFIP{a}{c}$ is an MDFIP, by Theorem~\ref{thm:charMDFIPs} we need to show that
$c \wedge a \prec a$ and $c \prec c \vee a$. By (JM-LSM) we only need to prove $c \prec c \vee a$, which follows from Lemma~\ref{lem:Tab}.
We have shown that (L-abc) holds.

Now assume
(L-abc). To show (JM-LSM), let $a \in \mathsf{J}(L)$, $b \in \mathsf{M}(L)$ and $b \prec a \vee b$. We need to prove
$a \wedge b \prec a$. From $b \prec a \vee b$ we have $a \nleqslant b$. By (L-abc) there exists $c \geqslant b$ such that $\MDFIP{a}{c}$ is an MDFIP. Hence
$c \in \mathsf{M}(L)$ and by Theorem~\ref{thm:charMDFIPs},
$c\wedge a \prec a$.
We claim that $c = b$.
Suppose
that $c>b$. Then, since $b \in \mathsf{M}(L)$, it has a unique upper cover $b^{\star}$. As $b \prec a \vee b$, we get $b^{\star} = a \vee b$.
From $c>b$
we have $c \ge b^{\star}= a \vee b \ge a$. This contradicts the fact that $\MDFIP{a}{c}$ is an MDFIP. Hence
$c=b$.
This proves $a \wedge b = c\wedge a \prec a$
as required.
\end{proof}

\begin{rem}
We notice that if a finite lattice $L$ satisfies
(L-abc), then in the situation $a \nleqslant b$ for $a \in \mathsf{J}(L)$, $b \in \mathsf{M}(L)$,
an arbitrary maximal element of $T_{ab}$
can be taken
for the element $c \geqslant b$ such that
$\MDFIP{a}{c}$ is an MDFIP.
Indeed, if $c$ is any maximal element of $T_{ab}$, then $c \in \mathsf{M}(L)$, $a \nleqslant c$, $b \leqslant c$ and so by
the assumed condition (L-abc) there is
$c' \geqslant c$
such that $\MDFIP{a}{c'}$ is an MDFIP. Hence $c'\in \mathsf{M}(L)$, $a \nleqslant c'$, $b \leqslant c'$, thus $c' \in T_{ab}$. From the maximality of $c$ in $T_{ab}$ we get $c=c'$ as required.
\end{rem}

Now
we present a digraph condition dual to (JM-LSM). The condition below is a
strengthening of the (Ti) condition, and
because of its connection to lower semimodularity, we have chosen the name (LTi). Later, in Definition~\ref{def:UTi},  (UTi) is used for the dual condition related to upper semimodularity.
\begin{df}\label{def:LTi}
Consider the following condition on a
TiRS digraph $G=(V,E)$:
$$\text{\upshape (LTi)} \qquad uEv \:\Longrightarrow\: (\exists\, w \in V) (wE = uE \:\&\: Ew \subseteq Ev ).$$
\end{df}

Note that (LTi) is not dual to (LSM) as
Figure~\ref{fig:L3dualdigraph} shows.
The next two results prove that it is (JM-LSM)
that is dual to (LTi).

\begin{prop}\label{prop:labc-iff-LTi}
A finite TiRS digraph satisfies {\upshape (LTi)} if and only if it is the dual digraph of a lattice that satisfies {\upshape (L-abc)}.
\end{prop}

\begin{proof}
Assume a finite lattice $L$ satisfies (L-abc). To show that the dual digraph $G_L$ satisfies {\upshape (LTi)}, let $u=\MDFIP{a}{m}$, $v=\MDFIP{j}{b}$
be vertices of the
digraph
$G$ and let $uEv$, whence  $a \nleqslant b$. Then by (L-abc) there exists $c \in \mathsf{M}(L)$ such that $b\le c$ and
$\MDFIP{a}{c}$ is an MDFIP. If we denote $w=\MDFIP{a}{c}$ as a vertex of
$G$, then
by Lemma~\ref{lem:xEyE}
we have $wE = uE$ and $Ew \subseteq Ev$ as required.

For the converse,
assume that a finite TiRS digraph $G$ satisfies {\upshape (LTi)}. To show that its dual lattice $L$ satisfies (L-abc), let $a \in \mathsf{J}(L)$, $b \in \mathsf{M}(L)$ and $a \nleqslant b$.
Since $a \in \mathsf{J}(L)$ and $L$ is finite, by Proposition~\ref{prop:ambj}(i), there exists an element $m \in M(L)$ such that $u=\MDFIP{a}{m}$ is an MDFIP.
Similarly, since $b \in \mathsf{M}(L)$, by Proposition~\ref{prop:ambj}(ii) there exists $j \in J(L)$ such that $v = \MDFIP{j}{b}$ is an MDFIP.
Since $a \nleqslant b$,
we have $uEv$. Now, by {\upshape (LTi)}, there is a vertex $w=\MDFIP{c}{d} \in V(G)$ satisfying $wE = uE$ and $Ew \subseteq Ev$. Since $wE = uE$,
we get ${\uparrow} c = {\uparrow} a$,
so $c=a$. Since $Ew \subseteq Ev$,
Lemma~\ref{lem:xEyE}(ii) tells us that $d \geqslant b$.
This proves that $d$ is the desired element such that  $\MDFIP{a}{d}$ is an MDFIP.
\end{proof}

The main theorem of this section follows directly from Theorem~\ref{thm:labc-iff-JMLSM} and Proposition~\ref{prop:labc-iff-LTi}.

\begin{thm}\label{thm:LTi-iff-JM-LSM}
A finite TiRS digraph
is the dual digraph of a
finite lattice satisfying {\upshape (JM-LSM)} if and only if it
satisfies {\upshape (LTi)}.
\end{thm}

For completeness, we now state the conditions and results related to finite upper semimodular lattices and their dual digraphs.

\begin{df}
Let $L$ be a finite lattice. We say that $L$ satisfies the condition (JM-USM) if whenever $a \in \mathsf{J}(L)$, $b \in \mathsf{M}(L)$, and $a \wedge b \prec a$, then $b \prec a \vee b$. We say that $L$ satisfies
(U-abc) if whenever $a \in \mathsf{J}(L)$ and $b \in \mathsf{M}(L)$ and $a \nleqslant b$ then there exists $c \le a$ such that $\MDFIP{c}{b}$ is an MDFIP.
\end{df}

The proposition below connects the two conditions defined above.
\begin{prop}
A finite lattice satisfies {\upshape (U-abc)} iff it satisfies {\upshape (JM-USM)}.
\end{prop}

Our last definition is the condition (UTi) which is, like (LTi), a strengthening of the (Ti) condition from Definition~\ref{def:TiRS}.
\begin{df}\label{def:UTi}
Consider the following condition on a finite TiRS digraph $G=(V,E)$:
$$\text{\upshape (UTi)} \qquad uEv \:\Longrightarrow\: (\exists\, w \in V) (wE \subseteq uE \:\&\: Ev=Ew ).$$
\end{df}

\begin{thm}
A finite TiRS digraph satisfies {\upshape (UTi)} if and only if it is the dual digraph of a finite lattice that satisfies
{\upshape (JM-USM)}.
\end{thm}

\section{Dual digraphs of meet-distributive lattices}\label{sec:meet-dist}

In this section we will combine the results from Section~\ref{sec:lowsemimod} with results about dual digraphs of finite join- and meet-semidistributive lattices from~\cite{P5}.
The goal is give a description of the dual digraphs of finite meet-distributive lattices. This will give a description of a new class of structures that are in a one-to-one correspondence with the class of finite convex geometries. First, we recall some basic definitions.

A lattice $L$ is \emph{join-semidistributive} if it satisfies the following quasi-equation for all $a,b,c \in L$:
$$
\text{(JSD)} \qquad
a\vee b \approx a\vee c \quad \longrightarrow \quad a \vee b \approx a \vee (b \wedge  c).
$$
A lattice $L$ is \emph{meet-semidistributive} if it satisfies the following quasi-equation for all $a,b,c \in L$:
$$
\text{(MSD)} \qquad
a\wedge  b \approx a\wedge c \quad \longrightarrow \quad a \wedge b \approx  a \wedge (b \vee c).
$$
A lattice is \emph{semidistributive} if it satisfies both (JSD) and (MSD).

Considering the lattices in Figure~\ref{fig:dual-examples} one can see that  $N_5$ is semidistributive, $L_4$ is meet-semidistributive but not join-semidistributive, and $L_4^\partial$ is join-semidistributive but not meet-semidistributive.

For a finite lattice $L$ and $a \in L$, consider
$\mu(a) =\bigwedge  \{\,b \in L \mid b \prec a \,\}$.
A finite lattice is \emph{meet-distributive} (also called \emph{locally distributive}) if for any $a \in L$, the
interval $[\mu(a),a]$ is a distributive lattice
(c.f.~\cite[Section 5-2]{AN-Ch5}).

The following equivalence is extracted from \cite[Theorem 5-2.1]{AN-Ch5}.

\begin{thm}\label{thm:MD=JSD+LSM}
Let $L$ be a finite lattice. Then the following are equivalent:
\begin{enumerate}[label=\normalfont(\roman*)]
\item $L$ is meet-distributive;
\item $L$ satisfies {\upshape (JSD)} and {\upshape (LSM)}.
\end{enumerate}
\end{thm}

The
results below
use Theorem~\ref{thm:MD=JSD+LSM} to provide an additional characterisation of meet-distributive lattices using (JM-LSM), the condition that was central to Section~\ref{sec:lowsemimod}.
Later, we will use this to characterise their dual digraphs.

\begin{thm} \label{JM-LSM+JSD->LSM}
If a finite lattice $L$ satisfies {\upshape (JM-LSM)} and {\upshape (JSD)}, then it is lower semimodular.
\end{thm}

\begin{proof} Let $L$ be a finite lattice satisfying (JM-LSM) and (JSD). Let $a,b \in L$ be arbitrary such that $a \prec a \vee b$. We are going to show that $a \wedge b \prec b$. We will proceed by contradiction.

Suppose that $a \wedge b \nprec b$.
Since
$L$ is finite, there exists $c \in L$ such that
$a \wedge b < c < b$.
Then $b \nleqslant c$ and by Lemma~\ref{lem:JM} the set $S_{cb} = \{\, j \in \mathsf{J}(L) \mid j \leqslant b, j \nleqslant c\,\}$ is non-empty. Let $p$ be a minimal element of $S_{cb}$.

Suppose $p \leqslant a$, then since $p \le b$, we get $p \le a \wedge b \le c$, which is a contradiction, so $p \nleqslant a$. Then by Lemma~\ref{lem:JM}, the set $T_{pa} = \{ m \in M(L)$ $|$ $a \le m$ and $p \nleqslant m \}$ is non-empty. Let $m$ be a maximal element of $T_{pa}$. By Lemma~\ref{lem:Tab}, $m \prec m \vee p$. Since $m \in M(L)$, $p \in J(L)$, and $L$ satisifies (JM-LSM), we obtain $m \wedge p \prec p$.

It is easy to see that $ c \wedge p  < p$ as if  $c \wedge p = p$, then $p \le c$, which is a contradiction.

We will show in several steps that $m \wedge p < c \wedge p$. Firstly, we will show $m \wedge p \leqslant  c \wedge p$. Suppose $m \wedge p \nleqslant c$.
By Lemma~\ref{lem:JM}
there exists $j \in J(L)$ satisfying $j \le m \wedge p$ and $j \nleqslant c$. Then $ j \le  p \le b$, so $j \le b$. Since $p$ is a minimal element of $S_{cb}$,
and $j$ is also in $S_{cb}$, we obtain $p=j$. Then $ p \le m \wedge p$, so $p = m \wedge p$. Hence $p \le m$, which is a contradiction. Therefore $m \wedge p \le c$. Since $m \wedge p \le p$, we have $m \wedge p \le c \wedge p$.

To show $m \wedge p < c \wedge p$,
suppose to the contrary that $m \wedge p = c \wedge p$. 
We will continue by showing that $a \vee c = a \vee b = a \vee p$.

Since $p \le b$ we have $a \le a \vee p \le a \vee b$, and since $a \prec a \vee b$, we get $a = a \vee p$ or $a \vee p = a \vee b$. But $a \neq a \vee p$ since $p \nleqslant a$, so $a \vee b = a \vee p$.
Also, since $c \le b$, we have $a \le a \vee c \le a \vee b$. If $a \vee c = a$, then $c \le a$, whence $c \le a \wedge b$, which contradicts $a \wedge b < c$.
So $a \vee c = a \vee b = a \vee p$.
Hence, $m \vee c = (m \vee a) \vee c = m \vee (a \vee c) = m \vee (a \vee p) = (m \vee a) \vee p = m \vee p$. 

Now, by (JSD), $m \vee c = m \vee p = m \vee (c \wedge p) = m \vee (m \wedge p) = m$, which contradicts $p \nleqslant m$. This shows that
$m \wedge p < c \wedge p$ as required.

Now we have $m \wedge p < c \wedge p < p$. This  contradicts $m \wedge p \prec p$. Hence the element $c$ cannot exist, which shows that $a \wedge b \prec b$.
\end{proof}

\begin{rem}\label{rem:weakJSD}
Notice in the proof we actually use a weaker form of
(JSD).
We will say that a lattice $L$ is \emph{weakly join-semidistributive} if it satisfies the following quasi-equation for all $a \in \mathsf{M}(L)$, $b\in\mathsf{J}(L)$, $c \in L$:
$$
\text{(W-JSD)} \qquad
a\vee b \approx a\vee c \quad \longrightarrow \quad a \vee b \approx a \vee (b \wedge  c).
$$
Hence in Theorem~\ref{JM-LSM+JSD->LSM} we actually showed that (JM-LSM) and (W-JSD) implies (LSM).

We notice the lattice in Figure~\ref{fig:L3dualdigraph} satisfies (JM-LSM) but not (W-JSD): indeed $c \in \mathsf{M}(L)$, $b\in \mathsf{J}(L)$ and $c \vee b = c \vee a$ but $c \vee (b \wedge a) \neq c \vee a$.
\end{rem}

The result  below follows from
Theorems~\ref{thm:MD=JSD+LSM}
and~\ref{JM-LSM+JSD->LSM}.

\begin{coro}\label{coro:MD=JM-LSM+JSD}
A finite lattice is meet-distributive if and only if it satisfies both {\upshape (JM-LSM)} and {\upshape (JSD)}.
\end{coro}


The following theorem
provides a characterisation of the dual digraphs  of join- and meet-semidistributive lattices.
Notice that each of the conditions (i), (ii) and (iii) below is a strengthening of the (S) condition from the definition of TiRS digraphs (Definition~\ref{def:TiRS}).
\begin{thm}[{\cite[Theorem 3.6]{P5}}]\label{SD:char}
Let $G=(V,E)$
be a finite TiRS digraph with $u,v \in V$. Then
\begin{enumerate}[label=\normalfont(\roman*)]
\item $G$ is the dual digraph of a finite lattice satisfying {\upshape (JSD)} if and only if it satisfies the following condition:
$$ \text{\upshape (dJSD)} \qquad \text{if } u \neq v \text{ then } Eu \neq Ev.
$$
\item $G$ is the dual digraph of a finite lattice satisfying {\upshape (MSD)} if and only if it satisfies the following condition:
$$\text{\upshape (dMSD)} \qquad \text{if } u \neq v \text{ then } uE \neq vE. $$
\item $G$ is the dual digraph of a finite semidistributive lattice if and only if it satisfies the following condition:
$$\text{\upshape (dSD)} \qquad  \text{if } u \neq v \text{ then } Eu\neq Ev \text{ and } uE \neq vE. $$
\end{enumerate}
\end{thm}

The next few results in this section link the properties discussed earlier  to distributivity in lattices and transitivity in dual digraphs.

\begin{thm}\label{thm:LTi-to-trans}
Let $G=(V,E)$ be a finite TiRS digraph that satisfies both
{\upshape (dMSD)}
and {\upshape (LTi)}. Then $E$ is transitive.
\end{thm}

\begin{proof}
We first claim that if a finite TiRS digraph $G=(V,E)$ satisfies both (dMSD) and (LTi), then for any vertices $u,v\in V$, $uEv$ implies $Eu \subseteq Ev$. Indeed, $uEv$ by (LTi) implies the existence of $w\in V$ such that $wE = uE$ and $Ew \subseteq Ev$. By the property (dMSD),  $wE = uE$ means $w=u$, whence $Eu \subseteq Ev$ as required.

Now to show the transitivity of $E$, if $uEv$ and $vEw$ for some vertices $u,v,w\in V$, then by the above claim, $Eu \subseteq Ev$ and $Ev \subseteq Ew$. Hence $Eu \subseteq Ew$,  which means $u\in Ew$, whence $uEw$ as required.
\end{proof}

\begin{prop}\label{prop:trans-implies-poset}
If $G=(V,E)$ is TiRS digraph with transitive $E$,   then $G$ is a poset.
\end{prop}

\begin{proof}
As in a TiRS digraph $G=(V,E)$ the relation $E$ is reflexive, it only remains to show the antisymmetry of $E$.

Assume for $x,y \in V$ that $xEy$ and $yEx$. We firstly show that $xE \subseteq yE$: if $z\in V$ and $z\in xE$, then $xEz$ and with $yEx$  we get $yEz$ by transitivity of $E$, hence $z\in yE$ as required. Now $xE \subset yE$ by the condition (R) would give $(x,y)\notin E$, a contradiction. Hence $xE = yE$.

Analogously one can show that $Ey \subseteq Ex$ and since $Ey \subset Ex$ would by the condition (R) give $(x,y)\notin E$, we have $Ey = Ex$.
Using that $G$ satisfies the separation property (S), it follows that $x = y$ as required.
\end{proof}

The result below follows from
Birkhoff's representation,
Theorem~\ref{thm:LTi-to-trans} and Proposition~\ref{prop:trans-implies-poset}.

\begin{coro}\label{MSD+LTi->Dist}
If $L$
satisfies {\upshape (MSD)} and
{\upshape (JM-LSM)}, then $L$ is distributive.
\end{coro}

We now return to focus on finite meet-distributive lattices, with the goal of describing a class of digraphs connected to finite convex geometries.

Using the TiRS conditions, our conditions for the dual digraphs of (JM-LSM) and (JSD), respectively, and Corollary \ref{coro:MD=JM-LSM+JSD}, we get the following dual condition for meet-distributivity. Notice how
(dJSD)
is a strengthening of the (S) condition, and (LTi) is
a strengthening of the (Ti) condition.

 \begin{thm}\label{thm:dualMD}
 A finite digraph $G=(V,E)$ with a reflexive relation $E$ is the dual digraph of some finite meet-distributive lattice if and only if $G$ satisfies the following conditions:
\begin{enumerate}[align=parleft,leftmargin=*,labelsep=1.5cm,label=\normalfont]
\item[{\normalfont(dJSD)}]  If $x,y \in V$ and $x\neq y$ then $Ex \neq Ey$.
\item[{\normalfont(R)}] For all $x,y \in V$, if $xE \subset yE$ then $(x,y)\notin E$, and if $Ey \subset Ex$ then $(x,y)\notin E$.
\item[{\normalfont(LTi)}] For all $x,y \in V$, if $xEy$ then there exists $z \in V$ such that $zE = xE$ and $Ez \subseteq Ey$.
\end{enumerate}
 \end{thm}

\begin{proof}
Let $G$ be the dual digraph of some finite  meet-distributive lattice $L$. Then by Theorem~\ref{thm:TiRSrep} the digraph $G$ will satisfy (R). By Corollary~\ref{coro:MD=JM-LSM+JSD}, $L$ satisfies (JSD) and (JM-LSM). Hence by Theorem~\ref{SD:char}(i), $G$ satisfies (dJSD). Lastly, by Theorem~\ref{thm:LTi-iff-JM-LSM}, $G$ will satisfy (LTi).

Conversely, assume $G$ satisfies (dJSD), (R) and (LTi). Clearly $G$ is a TiRS digraph, hence the dual of a finite lattice $L$.
Theorem~\ref{SD:char}(i) shows that $L$ satisfies (JSD) and Theorem~\ref{thm:LTi-iff-JM-LSM} implies that $L$ satisfies (JM-LSM). Hence by Corollary~\ref{coro:MD=JM-LSM+JSD}, $L$ is meet-distributive.
\end{proof}

The theorem above establishes a one-to-one correspondence between finite meet-distributive lattices and finite digraphs satisfying the conditions (dJSD), (R) and (LTi). It is a restriction of Theorem~\ref{thm:TiRSrep}, while still generalising Birkhoff's
one-to-one correspondence between finite distributive lattices and finite posets.

\begin{df}[{\cite[Definition 30]{LTF}}]\label{def:closure}
Let $X$ be a set and $\phi : \powerset (X) \to \powerset (X)$. Then $\phi$ is a \emph{closure operator} on $X$ if for all $Y,Z \in \powerset(X)$
\begin{enumerate}[label=\normalfont(\roman*)]
\item $Y \subseteq \phi(Y)$;
\item $Y \subseteq Z$ implies $\phi(Y) \subseteq \phi(Z)$;
\item $\phi(\phi(Y))=\phi(Y)$.
\end{enumerate}
If $X$ is a set and $\phi$ a closure operator on $X$ then the pair $\langle X,\phi \rangle$ is called a \emph{closure system}.
For $Y\subseteq X$ we say that $Y$ is \emph{closed} if $\phi(Y)=Y$. The closed sets of a closure operator $\phi$ on $X$ form a complete lattice, denoted by $\Cld(X,\phi)$.
A \emph{zero-closure} system is a closure system $\langle X, \phi \rangle$ such that $\phi(\emptyset)=\emptyset$.
\end{df}

Now we turn our attention to convex geometries. The presentation here follows that of the book chapter by Adaricheva and Nation~\cite{AN-Ch5}.

\begin{df}[{\cite[Definition 5-1.1]{AN-Ch5}}]
A closure system $\langle X, \phi \rangle$ satisfies the \emph{anti-exchange property} if for all $x \neq y$ and all closed sets $A \subseteq X$,
$$\text{(AEP)} \qquad
x \in \phi (A \cup \{y\}) \text{ and } x \notin A \text{ imply that } y \notin \phi (A \cup \{x\}).
$$
\end{df}

\begin{df}[{\cite[Definition 1.6]{AGT-03}}]
A zero-closure system that satisfies the anti-exchange property is called a \emph{convex geometry}.
\end{df}

We now combine
Theorem~\ref{thm:dualMD} with
known equivalences
to obtain the following characterisation of finite convex geometries.
There are other equivalent
conditions~\cite[Theorem 5-2.1]{AN-Ch5}
that we have not included here.
\begin{thm}\label{thm:equivFCG} Let $L$ be a finite lattice. Then the following  are equivalent:
\begin{enumerate}[label=\normalfont(\roman*)]
\item $L$ is the closure lattice $\mathrm{Cld}(X,\phi)$ of a closure space $\langle X, \phi \rangle $ with the {\upshape (AEP)}.
\item $L$ is a meet-distributive lattice.
\item $L$ satisfies {\upshape(JSD)} and {\upshape(LSM)}.
\item $L$ satisfies {\upshape (JSD)} and {\upshape(JM-LSM)}.
\item $L$ is the lattice $\mathbb{C}(G)$ of a reflexive digraph $G$ satisfying {\upshape(dJSD)}, {\upshape(R)} and {\upshape (LTi)}.
\end{enumerate}

\end{thm}
\begin{proof} The equivalences of (i), (ii) and (iii)
are known~\cite[Theorem 5-2.1]{AN-Ch5}.
The equivalence of (iii) and (iv) is the result of Corollary~\ref{coro:MD=JM-LSM+JSD}, and the equivalence of (iv) and (v) is Theorem~\ref{thm:dualMD}.
\end{proof}

\section{Dual digraphs of finite modular lattices}\label{sec:mod} 

In this section we provide two sufficient conditions for a finite TiRS digraph to be the dual digraph of a finite modular lattice.

For $i=0,1,2$, let us denote by $G_i=(V_i,E_i)$ an induced subgraph of $G_{N_5}$
%
(see Figure~\ref{fig:dual-examples})
with $V_i=\{x,y,z\}$ and with $i$ of the arcs $xEy$ and $yEz$ missing compared to $G_{N_5}$. (For $i=1$ we can, w.l.o.g., consider the arc $yEz$ missing.) Hence $G_0 =G_{N_5}$, $G_1$ has one arc and an isolated vertex, and $G_2$ has no arc and consists of two isolated vertices. All three digraphs are reflexive, hence they have loops at each vertex.

We introduce the following condition for the dual digraph $G_L$ of a finite lattice $L$ in terms of ``Forbidden Induced Subgraphs'':
$$
\text{\upshape (FIS)}\qquad  \text{$G_L$ has neither $G_0=G_{N_5}$ nor $G_1$ as an induced subgraph.}
$$
The next lemma and two propositions lead to showing that the condition (FIS) is sufficient for modularity of a finite lattice $L$.
Note that by Lemma~\ref{lem:JM}, for $a, b \in L$ with $a \nleqslant b$, there always exist elements $\underline{a} \leqslant a$ and $\overline{b} \geqslant b$ such that $\MDFIP{\underline{a}}{\overline{b}}$ is an MDFIP. 
Below we 
write $a||b$ to indicate that $a\nleqslant b$ and $b \nleqslant a$.

\begin{lem}\label{lem:mod-sufficiency}
Let $a, b, c, 0, 1$ be any elements of the lattice that form a sublattice isomorphic to $N_5$
{\upshape(}where $0<a, b, c<1$, $b<c$ and $a||b, a||c${\upshape)}.
Let $x=\langle {\uparrow}\underline{a}, {\downarrow}\overline{c}\rangle$, $y=\langle {\uparrow}\underline{c}, {\downarrow}\overline{b}\rangle$ and $z=\langle {\uparrow}\underline{b}, {\downarrow}\overline{a}\rangle$ be any maximal disjoint extensions of $\langle {\uparrow}a, {\downarrow}c\rangle$, $\langle {\uparrow}c, {\downarrow}b\rangle$ and $\langle {\uparrow}b, {\downarrow}a\rangle$, respectively. Then the induced subgraph $\{x,y,z\}$ of $G_L$ is isomorphic either to $G_0 = G_{N_5}$, $G_1$, or $G_2$.
\end{lem}
\begin{proof}
First we must confirm that $x,y,z$ are distinct MDFIPs. If $x=y$ then ${\uparrow}\underline{a}={\uparrow}\underline{c}$ which implies ${\uparrow}\underline{a} \cap {\downarrow}\overline{c}\neq \emptyset$, i.e $x$ would not be an MDFIP. If $x=z$ then ${\uparrow}\underline{a}={\uparrow}{\underline{b}}$ which means $z$ would not be an MDFIP. Lastly, if $y=z$ then ${\downarrow}\overline{b}={\downarrow}\overline{a}$ and $z$ would not be an MDFIP. 

We claim that in the induced subgraph $\{x,y,z\}$ of $G_L$, the arcs  $xEy$ and $yEz$ are possible, but the induced subgraph $\{x,y,z\}$ has none of the other four possible arcs between distinct vertices: indeed, the arcs $yEx$, $zEy$, $xEz$ and $zEx$ are not present in $G_L$ because clearly
$c\in {\uparrow}\underline{c} \cap {\downarrow}\overline{c}$, $b\in {\uparrow}\underline{b} \cap {\downarrow}\overline{b}$, $a\in {\uparrow}\underline{a} \cap {\downarrow}\overline{a}$ and $b\in {\uparrow}\underline{b} \cap {\downarrow}\overline{c}$, respectively.

Hence $\{x,y,z\}$ is isomorphic to $G_i$ in case $i$ of the arcs $xEy$ and $yEz$ are missing in the induced subgraph $\{x,y,z\}$ for $i=0,1,2$.
\end{proof}

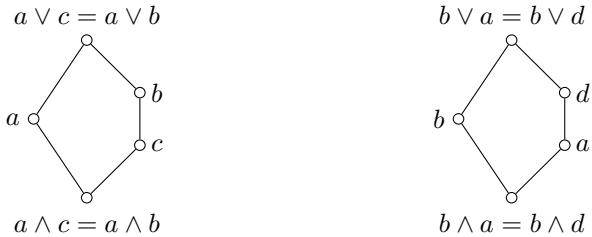
\begin{figure}[ht]
\centering
\begin{tikzpicture}[scale=0.7]
\begin{scope}[xshift=0cm]
\node[unshaded] (bot) at (0,0) {};
\node[unshaded] (a) at (-1,1.5) {};
\node[unshaded] (b) at (1,2) {};
\node[unshaded] (c) at (1,1) {};
\node[unshaded] (top) at (0,3) {};
\draw[order] (bot)--(a)--(top)--(b)--(c)--(bot);
\node[label,anchor=north] at (bot) {$a\wedge c=a\wedge b$};
\node[label,anchor=east,xshift=1pt] at (a) {$a$};
\node[label,anchor=west,xshift=-2pt] at (b){$b$};
\node[label,anchor=west,xshift=-2pt] at (c){$c$};
\node[label,anchor=south] at (top) {$a\vee c = a \vee b $};
\end{scope}

\begin{scope}[xshift=8cm]
\node[unshaded] (bot) at (0,0) {};
\node[unshaded] (a) at (-1,1.5) {};
\node[unshaded] (b) at (1,2) {};
\node[unshaded] (c) at (1,1) {};
\node[unshaded] (top) at (0,3) {};
\draw[order] (bot)--(a)--(top)--(b)--(c)--(bot);
\node[label,anchor=north] at (bot) {$b \wedge a = b \wedge d$};
\node[label,anchor=east,xshift=1pt] at (a) {$b$};
\node[label,anchor=west,xshift=-2pt] at (b){$d$};
\node[label,anchor=west,xshift=-2pt] at (c){$a$};
\node[label,anchor=south] at (top) {$b \vee a = b \vee d$};
\end{scope}
\end{tikzpicture}
\caption{The isomorphic copies of $N_5$ constructed in Proposition~\ref{prop:FIS->LSM} (left) and Proposition~\ref{prop:FIS->USM} (right).}\label{fig:FIS-N5}
\end{figure}

\begin{prop}\label{prop:FIS->LSM}
Let $L$ be a finite lattice and assume that its dual digraph $G_L = (V, E)$ satisfies {\upshape(FIS)}.
Then $L$ is lower semimodular.
\end{prop}
\begin{proof}
Suppose to the contrary that
$L$ does not satisfy (LSM).
Then there exist elements $a,b \in L$ such that $a \prec a \vee b$ but $a \wedge b \nprec b$.
Then
there exists an element $c \in L$ such that $a\wedge b <c<b$.
Hence 
$a\vee c \le a \vee b$. Since $a \prec a \vee b$, and $a \le a\vee c \le a \vee b$,
we get $a \vee c = a$ or $a \vee c = a \vee b$. If $a \vee c = a$, then $c \le a$,
so $c \le a \wedge b$, which contradicts $a \wedge  b < c$. It follows that $a \vee c = a \vee b$.
From $c < b$ we get $a \wedge c \leqslant a \wedge b$. Further, since $a \wedge b < c$ we get $a \wedge (a \wedge b)=a\wedge b \leqslant a \wedge c$. Thus $a \wedge c = a \wedge b$.

Hence
$a,c,b,a\wedge b, a\vee b$
forms a sublattice isomorphic to $N_5$ (see Figure~\ref{fig:FIS-N5}).
Let $x=\langle {\uparrow}\underline{a}, {\downarrow}\overline{b}\rangle$, $y=\langle {\uparrow}\underline{b}, {\downarrow}\overline{c}\rangle$ and $z=\langle {\uparrow}\underline{c}, {\downarrow}\overline{a}\rangle$, be
arbitrary maximal disjoint extensions of $\langle {\uparrow}a, {\downarrow}b\rangle$, $\langle {\uparrow}b, {\downarrow}c\rangle$ and $\langle {\uparrow}c, {\downarrow}a\rangle$, respectively.
Then by Lemma~\ref{lem:mod-sufficiency},
the induced subgraph $\{x,y,z\}$ of $G_L$ is isomorphic to $G_0 = G_{N_5}$, $G_1$, or $G_2$.
Using the assumption (FIS)
$\{x,y,z\}$ must be isomorphic to $G_2$.

In particular, it follows that
$G_L$ does not have the arc
$yEz$. Therefore $\underline{b} \le \overline{a}$. Suppose $a= \overline{a}$. Then
$\underline{b} \le a$,
so $\underline{b} \le a \wedge b$. This gives $\underline{b}
\le c \le \overline{c}$, which
contradicts
the fact that $y=\langle {\uparrow}\underline{b}, {\downarrow}\overline{c}\rangle$ is a disjoint filter-ideal pair.
Hence $a<\overline{a}$. Now either
$\overline{a}<a\vee b$
or $\overline{a} || a \vee b$, since if $\overline{a} \ge a \vee b >c \ge \underline{c}$ then $z=\langle {\uparrow}\underline{c}, {\downarrow}\overline{a}\rangle$ could not be a disjoint filter-ideal pair.

If $a < \overline{a} < a \vee b$, this contradicts $a \prec a \vee b$, so $\overline{a} || a \vee b$. If $\underline{b} > a$ then $b>\underline{b} > a$, which contradicts $a || b$.
If $\underline{b} \le a$, then
$\underline{b} \le a\wedge b \le c \le \overline{c}$, which contradicts
that $y=\langle {\uparrow}\underline{b}, {\downarrow}\overline{c}\rangle$ is a disjoint filter-ideal pair. This proves that $\underline{b} || a$. Since $\underline{b}\le b$, $a \vee \underline{b}  \le a \vee b$. If $a \vee \underline{b} = a \vee b$, then since $a<\overline{a}$ and $\underline{b} \le \overline{a}$, we get
$\overline{a} \ge a \vee \underline{b} = a \vee b$, which contradicts $\overline{a} || a \vee b$.
This establishes that
$a \vee \underline{b} < a \vee b$ and
$a < a \vee \underline{b}$ (since $\underline b || a$), which contradicts $a \prec a \vee b$.
Hence, our assumption that $L$ does not satisfy (LSM) leads to a contradiction. \end{proof}

Below we give the result dual to Proposition~\ref{prop:FIS->LSM}.
The proof is similar to the above argument, so we omit some of the details.

\begin{prop}\label{prop:FIS->USM}
Let $L$ be a finite lattice and assume that its dual digraph $G_L = (V, E)$ satisfies {\upshape(FIS)}.
Then $L$ is upper semimodular.
\end{prop}

\begin{proof}
Suppose $L$ does not satisfy (USM).
Then there are elements $a,b \in L$ such that $a \wedge b \prec b$ but
$a \nprec a \vee b$, i.e. there is $d \in L$ such that $a < d<a \vee b$.
Analogous to the proof of Proposition~\ref{prop:FIS->LSM}, it can be shown that the elements $b,a,d,a\wedge b, a \vee b$ form a sublattice isomorphic to $N_5$ (see Figure~\ref{fig:FIS-N5}).

Then by Lemma~\ref{lem:mod-sufficiency}, arbitrary maximal disjoint extensions of $\langle {\uparrow}b, {\downarrow}d\rangle$, $\langle {\uparrow}d, {\downarrow}a\rangle$ and $\langle {\uparrow}a, {\downarrow}b\rangle$, denoted by $x=\langle {\uparrow}\underline{b}, {\downarrow}\overline{d}\rangle$, $y=\langle {\uparrow}\underline{d}, {\downarrow}\overline{a}\rangle$ and $z=\langle {\uparrow}\underline{a}, {\downarrow}\overline{b}\rangle$, respectively, form an induced subgraph $\{x,y,z\}$ of $G_L$ that is isomorphic either to $G_0 = G_{N_5}$, $G_1$, or $G_2$.
Using (FIS),
$\{x,y,z\}$ is isomorphic to $G_2$.

In particular, it follows that
$G_L$ does not have the arc $xEy$. Hence,  $\underline{b} \le \overline{a}$.
We can then get $\underline{b}<b$  (as we got $a<\overline{a}$ in 
Proposition~\ref{prop:FIS->LSM}---see 
the left lattice in Figure~\ref{fig:FIS-N5}). Now either $a \wedge b< \underline{b}$
or $a \wedge b || \underline{b}$.

If $a \wedge b < \underline{b} < b$, this contradicts $a \wedge b \prec b$, so $\underline{b} || a \wedge b$. We can also show $b || \overline{a}$ (as we showed $\underline{b} || a$ in Proposition~\ref{prop:FIS->LSM}).

Since $a \le \overline{a}$, 
we get 
$a \wedge b  \le \overline{a} \wedge b$.
We can again establish that
$a \wedge b  < \overline{a} \wedge b$ and
$\overline{a} \wedge b < b$ (since $b || \overline{a}$), which contradicts $a \wedge b \prec b$.
Hence, our assumption that $L$ does not satisfy (USM) leads to a contradiction. \end{proof}

Now we can deduce that the condition (FIS) is a sufficient condition for modularity of a finite lattice.

\begin{thm}\rm{(}{\bf Sufficient condition for modularity}\rm{)}\label{thm:mod-sufficiency}
Let $L$ be a finite lattice with dual TiRS digraph
$G_L$.  If $G_L$ satisfies the condition
{\upshape(FIS)}
then $L$ is modular.
\end{thm}
\begin{proof}
If follows by Propositions~\ref{prop:FIS->LSM} and~\ref{prop:FIS->USM} that $L$ satisfies both (LSM) and (USM). Since $L$ is finite, we have that $L$ is modular.
\end{proof}

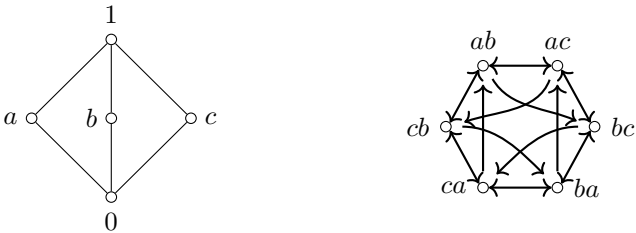
\begin{figure}[ht]
\centering
\begin{tikzpicture}[scale=0.7]
\begin{scope}[xshift=0cm]
\node[unshaded] (bot) at (0,0) {};
\node[unshaded] (a) at (-1.5,1.5) {};
\node[unshaded] (b) at (0,1.5) {};
\node[unshaded] (c) at (1.5,1.5) {};
\node[unshaded] (top) at (0,3) {};
\draw[order] (bot)--(a)--(top)--(b)--(bot)--(c)--(top);
\node[label,anchor=east,xshift=1pt] at (a) {$a$};
\node[label,anchor=east,xshift=1pt] at (b) {$b$};
\node[label,anchor=west,xshift=-1pt] at (c) {$c$};
\node[label,anchor=north] at (bot) {$0$};
\node[label,anchor=south] at (top) {$1$};
 \end{scope}

\begin{scope}[yshift=2.5cm,xshift=7cm]
\node[unshaded] (ab) at (0,0) {};
\node[label,anchor=south] at (ab) {$ab$};
\node[unshaded] (ac) at (1.4,0) {};
\node[label,anchor=south] at (ac) {$ac$};
\node[unshaded] (ba) at (1.4,-2.32) {};
\node[label,anchor=west] at (ba) {$ba$};

\node[unshaded] (bc) at (2.1,-1.16) {};
\node[label,anchor=west] at (bc) {$bc$};
\node[unshaded] (ca) at (0,-2.32) {};
\node[label,anchor=east] at (ca) {$ca$};
\node[unshaded] (cb) at (-0.7,-1.16) {};
\node[label,anchor=east] at (cb) {$cb$};

\path[thick,<->] (ab.east) edge  (ac.west);
\path[thick,<->] (ac.south east) edge (bc.north west);
\path[thick,<->] (bc.south west) edge (ba.north east);
\path[thick,<->] (ba.west) edge  (ca.east);
\path[thick,<->] (ca.north west) edge  (cb.south east);
\path[thick,<->] (cb.north east) edge  (ab.south west);
\path[thick,->,shorten >=5pt,shorten <=4pt] (ab.south east) [out=300,in=155] edge  (bc.west);

\path[thick,->,shorten >=5pt,shorten <=4pt] (bc.west) [out=180,in=45] edge (ca.north east);
\path[thick,->,shorten >=5pt,shorten <=4pt] (ca.north) edge (ab.south);
\path[thick,->,shorten >=5pt,shorten <=4pt] (ac.south west) [out=250,in=20] edge  (cb.east);
\path[thick,->,shorten >=5pt,shorten <=4pt] (cb.east) [out=0,in=135] edge  (ba.north west);
\path[thick,->,shorten >=5pt,shorten <=4pt] (ba.north) edge  (ac.south);
\end{scope}
\end{tikzpicture}
\caption{$M_3$ and its dual digraph.}\label{fig:M3digraph}
\end{figure}

We notice that the dual digraph of the modular lattice $M_3$ has neither $G_0=G_{N_5}$ nor $G_1$ as an induced subgraph (see Figure~\ref{fig:M3digraph}), hence it satisfies (FIS).
The following example shows that the digraphs $G_0$ and $G_1$
cannot be dropped as forbidden induced subgraphs in the condition (FIS) for the dual digraph $G_L$, which guarantees the modularity of a finite lattice $L$.

\begin{eg}
The dual digraph of $L_3^{\partial}$ in  Figure~\ref{fig:L3dualdigraph} contains $G_0$ as an induced subgraph, but not $G_1$. Hence
the lattice $L_3^{\partial}$
(in addition to $N_5$)
witnesses
that the digraph $G_0$ cannot be dropped from the condition (FIS).

The dual digraphs of the lattices $L_4$ and $L_4^\partial$ in Figure~\ref{fig:dual-examples}
do not contain $G_0$ as an induced subgraph but they both contain $G_1$ as an induced subgraph. Hence these two examples witness that the digraph $G_1$ cannot be dropped from the condition (FIS).
\end{eg}

Now we are going to show that the condition (FIS) is not necessary for modularity. Indeed, it is not the case that every lattice whose dual digraph has $G_0=G_{N_5}$ as an induced subgraph is a non-modular lattice. The next example gives a modular lattice whose dual digraph has $G_0$ as an induced
subgraph (but does not have $G_1$ as an induced subgraph).

\begin{eg}\rm{(}{\bf Condition (FIS) not necessary for modularity}\rm{)}\label{mod:necessity}
Figure~\ref{fig:fis-counter}
shows a modular lattice $K$ on the left, and its dual digraph on the right. The
induced
subgraph isomorphic to $G_0$
is shown with the dotted arrows ($dcEcb$ and $cbEed$).
\begin{figure}[ht]
\centering
\begin{tikzpicture}[scale=0.7]
\begin{scope}[yshift=-5cm,xshift=4cm]
\node[unshaded] (bot) at (0,0) {};
\node[unshaded] (a) at (-1,2) {};
\node[unshaded] (c) at (1,2) {};
\node[unshaded] (b) at (0,2) {};
\node[unshaded] (e) at (0,1) {};
\node[unshaded] (top) at (0,3) {};
\node[unshaded] (d) at (-1, 1) {};
\draw[order] (e)--(b)--(top)--(c)--(e)--(a)--(top);
\draw[order] (e)--(bot)--(d)--(a);
\node[label,anchor=north] at (bot) {$0$};
\node[label,anchor=east,xshift=1pt] at (a) {$a$};
\node[label,anchor=west,xshift=-1pt] at (c) {$c$};
\node[label,anchor=west,xshift=-2pt] at (b){$b$};
\node[label,anchor=east,xshift=1pt] at (d){$d$};
\node[label,anchor=west,xshift=1pt] at (e){$e$};
\node[label,anchor=south] at (top) {$1$};
\end{scope}

\begin{scope}[yshift=-5cm,xshift=9cm]
\node[unshaded] (ca) at (3,0) {};
\node[label,anchor=north] at (ca) {$ca$};
\node[unshaded] (cb) at (1,0) {};
\node[label,anchor=north] at (cb) {$cb$};
\node[unshaded] (db) at (0,1.5) {};
\node[label,anchor=east] at (db) {$db$};
\node[unshaded] (dc) at (1,3) {};
\node[label,anchor=south] at (dc) {$dc$};
\node[unshaded] (bc) at (3,3) {};
\node[label,anchor=south] at (bc) {$bc$};
\node[unshaded] (ba) at (4,1.5) {};
\node[label,anchor=south,xshift=3pt,yshift=-1pt] at (ba) {$ba$};
\node[unshaded] (ed) at (6,1.5) {};
\node[label,anchor=south] at (ed) {$ed$};
\path[thick,shorten <=2.5pt, shorten >=2.5pt,<->] (ca.west) edge  (cb.east);
\path[thick,<->,shorten >=2pt,shorten <=2pt] (ca.north) edge (ba.south);
\path[thick,dashed,->,shorten >=6.5pt,shorten <=2.5pt] (cb.north east) edge (ed.south west);
\path[thick, shorten >=3pt, shorten <=4pt, <->] (ba.north west) edge (bc.south east);
\path[thick,shorten >=3pt, shorten <=4pt,<-] (ed.south) edge  (ca.north east);
\path[thick, shorten <=2.5pt, shorten >=2.5pt, ->] (ba.east) edge (ed.west);
\path[thick, shorten >=3.5pt, shorten <=3pt, ->] (bc.south east) edge (ed.north west);
\path[thick, shorten >=3pt, shorten <=4pt, <->] (dc.east) edge (bc.west);
\path[thick, shorten >=3pt, shorten <=4pt, <->] (dc.south west) edge (db.north west);
\path[thick, shorten >=3pt, shorten <=4pt, <->] (cb.north west) edge (db.south east);
\path[thick, shorten <=2.5pt, shorten >=2.5pt, ->] (ba.west) edge (dc.south east);
\path[thick, shorten <=2.5pt, shorten >=2.5pt, ->] (cb.north) edge (ba.west);
\path[thick, shorten <=2.5pt, shorten >=2.5pt, ->] (db.east) edge (bc.south west);
\path[thick, shorten <=2.5pt, shorten >=2.5pt, ->] (ca.north west) edge (db.east);
\path[thick, dashed, shorten <=2.5pt, shorten >=2.5pt, ->] (dc.south) edge (cb.north);
\end{scope}

\end{tikzpicture}
\caption{A finite modular lattice $K$ whose dual digraph contains $G_0=G_{N_{5}}$ as an induced subgraph.}
\label{fig:fis-counter}
\end{figure}
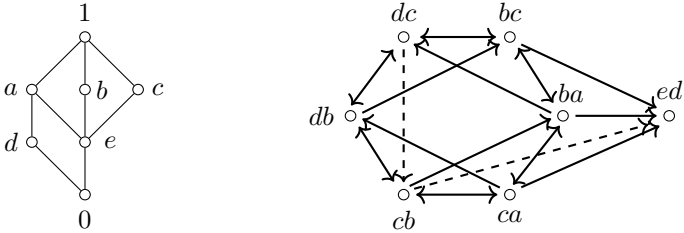

\end{eg}

The fact that the dual TiRS digraph $G_L=(V,E)$ of a finite modular lattice $L$ does not contain
$G_0=G_{N_5}$ as an induced subgraph can be
understood as some form of a ``weak transitivity'' condition for $G_L$. We cannot have the arcs
$xEy$ and $yEz$ in $G_L$ without having also
the arc  $xEz$ or at least the arc  $zEx$ (provided there are no ``opposite'' arcs  $yEx$ and $zEy$
in~$G_L$):
\begin{align*}
\text{(wT0)} \qquad &
\text{for all vertices $x,y,z\in V$, if $xEy$ and $yEz$, but $(y,x)\notin E$ and}\\
& \text{$(z,y)\notin E$, then $xEz$ or $zEx$.}
\end{align*}

Similarly, the fact that the dual TiRS digraph $G_L=(V,E)$ of a finite modular lattice $L$ does not contain the digraph $G_1$ as an induced subgraph can be
understood as some form of a ``weak transitivity'' condition for $G_L$:
\begin{align*}
\text{(wT1)} \qquad &
\text{for all vertices $x,y,z\in V$, if $xEy$ but $(y,x)\notin E$ and $(y,z)\notin E$}\\
&\text{and
$(z,y) \notin E$
then $xEz$ or $zEx$.}
\end{align*}

\begin{eg}
It is easy to see that the dual digraph of the lattice $M_3$
(Figure~\ref{fig:M3digraph})
satisfies the weak transitivity conditions (wT0) and (wT1).
The
lattices $L_4$ and $L_4^{\partial}$ in Figure~\ref{fig:dual-examples}, and $L_3^{\partial}$ in Figure~\ref{fig:L3dualdigraph}
are non-modular lattices.
The weak transitivity condition (wT0)
is not satisfied  in the dual digraph of
$L_3^{\partial}$. In the dual digraphs of the lattices  $L_4$ and $L_4^{\partial}$ we see the failures of (wT1).
\end{eg}

We notice that the weak transitivity conditions (wT0) and (wT1) are essentially expressing on the digraph side that the digraph $G_L$ does not contain respectively the graphs $G_0$ and $G_1$ as induced subgraphs.

Hence the sufficiency of the quasi-equations (wT0) and (wT1) on the dual TiRS digraphs $G_L$ for the modularity of $L$ comes
as no surprise:

\begin{coro}\label{cor:suff-WTi}\rm{(}{\bf Sufficient condition for modularity by ``weak transitivity''}\rm{)}\label{SC:weak-trans}
Let $L$ be a finite lattice with dual TiRS digraph $G_L = ( V, E )$. If $G_L$ satisfies the weak transitivity conditions (wT0) and (wT1), then $L$ is modular.
\end{coro}

\begin{proof}
Let the weak transitivity conditions (wT0) and (wT1) be satisfied in $G_L$. Suppose for contradiction that the lattice $L$ is not
modular. Then by Theorem~\ref{thm:mod-sufficiency}, for some $i\in\{0,1\}$ the digraph $G_L$ contains the digraph $G_i$ as an induced subgraph on
certain vertices $x,y,z\in V$. It follows that the  weak transitivity condition (wTi) is not satisfied.
\end{proof}

\section{Conclusions and future work}\label{sec:conclusion}

In this paper, we firstly (in Section~\ref{sec:lowsemimod})
defined
two
lattice
conditions which generalise lower semimodularity and (upper) semimodularity
respectively.
We were motivated by
Figure~\ref{fig:FCAbookdiagram},
taken from Ganter and Wille's book~\cite{FCA99}
(see also the PhD thesis of Reppe~\cite[Chapter 3.7]{Reppe-thesis}).
There, weakenings of (LSM) and (USM) are given using complicated conditions on standard contexts.   Our lattice-theoretic conditions on finite lattices that are weakenings of (LSM) and (USM), which we call (JM-LSM) and (JM-USM), seem to be simpler than the mentioned conditions in Figure~\ref{fig:FCAbookdiagram} and they are easily seen to be natural generalisations of (LSM) and (USM). Our focus was the generalisation of lower
semimodularity,
and we characterised the dual of (JM-LSM) on the dual digraphs of finite lattices. We think
the answer to the question below will be affirmative, yet investigating it is beyond the scope of this paper.

\begin{prob}
Are the top left and top right conditions in Figure~\ref{fig:FCAbookdiagram}, in terms of Formal Concept
Analysis~\cite{FCA99},
equivalent to {\upshape (JM-LSM)} and {\upshape (JM-USM)}
respectively?
\end{prob}

In Section~\ref{sec:meet-dist} we used the results of Section~\ref{sec:lowsemimod} to obtain a new characterisation of meet-distributive lattices in Theorem~\ref{thm:MD=JSD+LSM}. Combining this with
previous results~\cite{P5},
we
obtained a characterisation of the dual digraphs of finite meet-distributive lattices.
Theorem~\ref{thm:equivFCG} shows that we have identified a new class of structures that is in a one-to-one correspondence with finite convex geometries.

In Remark~\ref{rem:weakJSD} we gave a condition, (W-JSD), which is a weakening of join-semidistributivity.
The lattice $M_3$ satisfies (LSM) but not (W-JSD)
and hence shows that (LSM) is not equivalent to (JM-LSM) and (W-JSD).
This leads us to ask the following question.
\begin{prob}
Is
there another weakening of  {\upshape (JSD)} such that when it is combined with {\upshape (JM-LSM)}, this will be equivalent to {\upshape (LSM)}?
\end{prob}

Theorem~\ref{thm:dualMD} gave three conditions ((dJSD), (R) and (LTi)) on reflexive digraphs,  which characterise the dual digraphs of finite meet-distributive lattices. This leads to the posing of the following open problem.

\begin{prob}
Can the conditions {\upshape (dJSD)}, {\upshape (R)} and {\upshape (LTi)} be combined to give fewer, and possibly simpler, conditions?
\end{prob}

In Section~\ref{sec:mod} we introduced the condition (FIS) on dual digraphs and showed that it implies both lower and upper semimodularity of a finite lattice. Hence (FIS) was shown to be a sufficient condition for modularity of a finite lattice (Theorem~\ref{thm:mod-sufficiency}).
We also formulated a sufficient condition for modularity in different terms in Corollary~\ref{cor:suff-WTi}. The condition (FIS) was shown not to be necessary for modularity of a finite lattice and hence we raise the following open question.

\begin{prob} Is it possible to find forbidden induced subgraphs that characterise the dual digraphs of finite modular lattices in
an analogous
way
to how
$N_5$ characterises modularity?
\end{prob}

The task of representing structures (in our case digraphs) dual to finite modular lattices has proved to be very challenging. We note that in the setting of formal
contexts
dual to finite lattices, a condition dual to semimodularity has been obtained~(c.f.~item (4) of \cite[Theorem 42]{FCA99}). We have attempted to translate this condition to TiRS digraphs and the result was
a
complicated and
opaque
condition. We do not believe that the translation of this condition and its dual will yield a useful characterisation of the TiRS digraphs dual to finite modular lattices.

\subsection*{Acknowledgements}
The first author acknowledges
the hospitality of Matej Bel University
during a visit in August--September 2022, as well as the
National Research Foundation (NRF) of South Africa (grant
127266). The second author acknowledges his appointment as a Visiting Professor at the University of Johannesburg from June 2020
and the hospitality shown during a visit in July--August 2023. He further acknowledges support by Slovak VEGA
grant 1/0152/22. The authors would like to thank Jos\'{e} S\~{a}o Jo\~{a}o for useful discussions on these topics.


\end{document}